\documentstyle[12pt,amsmath,amssymb]{article}

\begin{document}
\newtheorem{lem}{Lemma}
\newtheorem{th}{Theorem}
\newtheorem{prop}{Proposition}
\newtheorem{rem}{Remark}
\newtheorem{define}{Definition}
\newtheorem{cor}{Corollary}

\allowdisplaybreaks

\newcommand{\ddGamma}{\overset{{.}{.}}{\Gamma}_X}
\newcommand{\ddG}{\overset{{.}{.}}{\Gamma}_{X,0}}
\newcommand{\D}{{\cal D}}
\newcommand{\N}{{\Bbb N}}
\newcommand{\C}{{\Bbb C}}
\newcommand{\Z}{{\Bbb Z}}
\newcommand{\R}{{\Bbb R}}
\newcommand{\Rp}{{\R_+}}
\newcommand{\eps}{\varepsilon}
\newcommand{\FFF}{\widehat{{\cal F}}_{\mathrm fin}(\D)}
\newcommand{\DFFF}{\widehat{{\cal F}}{}^*_{\mathrm fin}(\D)}
\newcommand{\om}{\omega}

\newcommand{\F}{{\cal F}}
\newcommand{\FF}{\F_{\mathrm fin}(\D)}
\newcommand{\DFF}{\F^*_{\mathrm fin}(\D)}
\newcommand{\Gamman}{\Gamma_{X,0}}
\newcommand{\supp}{\operatorname{supp}}
\newcommand{\la}{\langle}
\newcommand{\ra}{\rangle}
\newcommand{\const}{\operatorname{const}}

\newcommand{\ho}{\widehat\otimes}
\newcommand{\ot}{\otimes}
\newcommand{\fii}{\varphi}

\renewcommand{\emptyset}{\varnothing}
\renewcommand{\tilde}{\widetilde}
\newcommand{\rom}[1]{{\rm #1}}

\newcommand{\di}{\partial}
\renewcommand{\div}{\operatorname{div}}

\newcommand{\wom}[1]{{:}\,\om^{\ot #1}{:}}
\newcommand{\gom}[1]{{:}\,\gamma^{\ot #1}{:}}

\thispagestyle{empty}

$\mbox{}$\vspace{0mm}
\begin{center}
\LARGE\bf On a Spectral Representation \\[2mm]
\LARGE\bf for Correlation Measures\\[2mm]
in Configuration Space Analysis\\[20mm]
\large\bf Yuri M. Berezansky\\[3mm]
\large Inst.\ Math., NASU, 252601 Kiev, Ukraine\\[15mm]
\large\bf Yuri G. Kondratiev\\[3mm]
\large Inst.\ f.\ Angew.\ Math., Univ.\ Bonn, 53115 Bonn, Germany;\\
\large BiBoS, Univ.\ Bielefeld, 33615 Bielefeld, Germany; and\\
Inst.\ Math., NASU, 252601 Kiev, Ukraine\\[15mm]
\large\bf Tobias Kuna\\[3mm]
\large Inst.\ f.\ Angew.\ Math., Univ.\ Bonn, 53115 Bonn, Germany; and\\
\large BiBoS, Univ.\ Bielefeld, 33615 Bielefeld, Germany\\[15mm]
\large\bf Eugene  Lytvynov\\[3mm]
\large Inst.\ f.\ Angew.\ Math., Univ.\ Bonn, 53115 Bonn, Germany; and\\
\large BiBoS, Univ.\ Bielefeld, 33615 Bielefeld, Germany
\vspace{30mm}
\end{center}

\setcounter{page}{0}

\newpage
\thispagestyle{empty}

\begin{center}\bf Abstract\end{center}
\noindent\begin{small}

The paper is devoted to the study of configuration space analysis by using the projective spectral
theorem. For a manifold $X$, let $\Gamma_X$, resp.\ $\Gamma_{X,0}$ denote the space of all, resp.\ finite
configurations in $X$. The so-called $K$-transform, introduced by A.~Lenard, maps functions on $\Gamma_{X,0}$
into functions on $\Gamma_{X}$ and its adjoint $K^*$ maps probability measures on $\Gamma_X$ into $\sigma$-finite
measures on $\Gamma_{X,0}$. For a probability measure $\mu$ on $\Gamma_X$, $\rho_\mu:=K^*\mu$ is called the correlation
measure of $\mu$. We consider the inverse problem of existence of a probability measure $\mu$ whose correlation measure
$\rho_\mu$ is equal to  a given measure $\rho$. We introduce  an operation
of
$\star$-convolution of two functions on $\Gamma_{X,0}$ and suppose that the measure $\rho$ is $\star$-positive definite,
which enables us to introduce the Hilbert space
${\cal H}_\rho$ of functions on $\Gamma_{X,0}$ with the scalar product $(G^{(1)},G^{(2)})_{{\cal H}_{\rho}}=
\int_{\Gamma_{X,0}}(G^{(1)}\star\overline G{}^{(2)})(\eta)\,\rho(d\eta)$. Under a condition on the growth
of the   measure $\rho$ on the $n$-point configuration spaces, we construct the Fourier transform
in generalized joint eigenvectors of some special family $A=(A_\fii)_{\fii\in\D}$, $\D:=C_0^\infty(X)$, of commuting
selfadjoint operators in ${\cal H}_\rho$. We show that this
 Fourier transform is  a unitary between ${\cal H}_{\rho}$ and the
$L^2$-space $L^2(\Gamma_X,d\mu)$, where $\mu$ is the spectral
measure of $A$. Moreover, this unitary coincides with the
$K$-transform, while the measure $\rho$ is the correlation measure
of $\mu$. \vspace{5mm}

\noindent 2000 {\it AMS  Mathematics Subject Classification}.
Primary 60G57, 47A75. Secondary 60K35
\end{small}

\newpage

\newpage

{\it To the memory of Professor Yuri L. Daletsky}

\section{Introduction}
The configuration space $\Gamma_X$ over a (non-compact) Riemannian
manifold $X$ is defined as the set of all locally finite subsets (configurations)
in $X$. Such spaces as well as probability measures on them appear naturally
in several topics of mathematics and physics. Let us mention only the theory of point
processes \cite{KMM, DVJ}, classical statistical mechanics \cite{Preston, Georgii},
and nonrelativistic quantum field theory, e.g., \cite{Menia, MS75} and references therein.

An important tool in the study of configuration space analysis is the so-called $K$-transform.
Roughly speaking, this transform maps functions defined on the space $\Gamma_{X,0}$ of
finite configurations in $X$ into functions defined on the space $\Gamma_X$ of all configurations.
Interpreting the algebra of functions on $\Gamma_X$ as observables of our system, we may consider functions on
$\Gamma_{X,0}$ as quasi-observables, from which we are able to reconstruct observables by using the $K$-transform.
This special kind of observables is known in physics and called additive type observables, see \cite{Bog46}.
A.~Lenard was the first to recognize the operator nature of the $K$-transform \cite{Len1, Len2,
Len3}. Recently,
this theory was reanalized and further developed
in \cite{Kuna1, Kuna11, Kuna2, Kuna3}, where the reader can find also many applications
of this transform.

The adjoint $K^*$ of the $K$-transform, defined  by the formula
$$\int_{\Gamma_X} (KG)(\gamma)\,\mu(d\gamma)=\int_{\Gamma_{X,0}}G(\eta)(K^*\mu)(d\eta),$$
maps probability measures on $\Gamma_X$ into $\sigma$-finite measures on $\Gamma_{X,0}$, and
$\rho_\mu:=K^*\mu$ is called the correlation measure of $\mu$.

In several applications, a $\sigma$-finite measure $\rho$ on $\Gamma_{X,0}$
appears as a given object
and the problem is to show that this $\rho$ can be seen as
a correlation measure for a probability measure on
$\Gamma_X$.
Different types
of sufficient conditions were given for this to hold.
A.~Lenard
 \cite{Len2,Len3} used essentially a positivity condition for the correlation measure,
which allowed him to construct a  linear positive functional and
apply a version of the Riesz--Krein extension theorem. His conditions were also necessary.
O.~Macchi
\cite{Macke} (see also \cite{DVJ}) needed an additional condition in order to get
an explicit construction of the measure on $\Gamma_X$.

The present paper is also devoted to this problem. As a first step,
we utilize the idea of \cite{Kuna1,Kuna3},
 introducing the so-called $\star$-convolution on a
space of functions on $\Gamma_{X,0}$ and demanding that $\rho$
be $\star$-positive definite, that is,
\begin{equation}\label{dft}
\int_{\Gamma_{X,0}}(G\star\overline G)(\eta)\,\rho(d\eta)\ge0.\end{equation}
Unlike the approach of \cite{Kuna1, Kuna3}, where the authors prove a
Bochner type theorem, we use a spectral approach. The condition \eqref{dft}
enables us to  introduce in Section~2 the $\star$-convolution
Hilbert space ${\cal H}_\rho$ of functions on $\Gamma_{X,0}$ with the scalar product
$$(G^{(1)},G^{(2)})_{{\cal H}_\rho}:=\int_{\Gamma_{X,0}}(G^{(1)}\star
\overline G{}^{(2)})(\eta)\,\rho(d\eta).$$

Next, we follow the general strategy of representation  of positive
definite kernels and functionals by using
the projective spectral theorem, see \cite{BK,new1, BeLi,Ly}.
We consider in the space ${\cal H}_{\rho}$ a family $(A_\varphi)_{\varphi\in\D}$
of Hermitian operators defined by the formula
$$(A_\varphi G)(\eta):=(\varphi\star G)(\eta)$$
on an appropriate domain. Here, $\D:=C_0^\infty(X)$
is the nuclear space
of all $C^\infty$ functions on $X$ with compact support, and each $\fii\in\D$ is identified with the function on
$\Gamma_{X,0}$ given as follows: $\fii(\eta):=\fii(x)$ if $\eta=\{x\}$ and $\fii(\eta):=0$ if the number of points in
$\eta\in\Gamma_{X,0}$ is not equal to one.

Under a rather weak condition on  the measure  $\rho$,
 we show that the operators $A_\fii$ are essentially selfadjoint  in
${\cal H}_\rho$ and their closures $A_\fii^\sim$ constitute a family of commuting selfadjoint operators in
${\cal H}_\rho$. Moreover, these operators are shown to satisfy the conditions of the projective spectral theorem,
and the Fourier transform in generalized joint eigenvectors of the family  $(A_\fii^\sim)_{\fii\in\D}$  gives
a unitary isomorphism between ${\cal H}_\rho$ and an $L^2$-space $L^2(\D',d\mu)$, where $\D'$ is the dual of $\D$
and $\mu$ is the spectral measure of the family $(A_\fii^\sim)_{\fii\in\D}$. Under this isomorphism, each
operator $A_\fii^\sim$ goes over into the operator of multiplication by the monomial $\la\cdot,\fii\ra$.
The corresponding Parseval inequality gives the required spectral representation
of the functional determined by the measure $\rho$.

In Section 3, following an idea in \cite{Kuna1}, we prove,
under  an additional, natural condition on $\rho$,
that the measure $\mu$ is concentrated actually on
$\Gamma_X$. Notice that the configuration space can be considered
as a subset of $\D'$ by identifying any configuration from $\Gamma_X$ with a sum
of
delta functions having support in the points of the configuration. Moreover,
the unitary constructed in Section~2
coincides now with the $K$-transform,
since $\rho=\rho_\mu$ is the correlation measure of $\mu$.

Finally, let us stress that  the spectral approach not only gives
an alternative way to find sufficient conditions for
a measure to be a correlation measure, but also  gives a new understanding of the $K$-transform as a {\it unitary\/}
operator between the Hilbert spaces ${\cal H}_\rho$ and $L^2(\Gamma_X,d\mu)$
which has the form of  Fourier transform.

\section{The $\star$-convolution Hilbert space and the corresponding
Fourier transform}

Let $X$ be a connected, oriented $C^\infty$ (non-compact) Riemannian manifold.
We denote by $\D$ the space $C_0^\infty(X)$ of all real-valued
infinite differentiable functions on $X$ with compact support. This
space can be naturally endowed with a topology of a nuclear space, see
e.g.\ \cite{Die}.

Let $\F_0(\D):=\C$ and $\F_n(\D):=\D_\C^{\ho n}$, $n\in\N$,
where $\D_\C$ denotes the complexification of the real space
$\D$ and $\ho$ stands for the symmetric tensor product. Notice that
$\F_n(\D)$ is the complexification of the space of all real-valued $C^\infty$  symmetric
functions on $X^n$
with compact support. Next, we define
$$\FF:=\bigoplus_{n=0}^\infty\F_n(\D)$$
to be the topological direct sum of the spaces $\F_n(\D)$, i.e., an arbitrary element
$G\in\FF$ is of the form $G=(G^{(0)},G^{(1)},\dots,G^{(n)},0,0,\dots)
$,
where $G^{(i)}\in\F_i(\D)$, and the convergence in
$\FF$ means the uniform finiteness and the coordinate-wise convergence. In
what follows, we will identify a $G^{(n)}\in\F_n(\D)$ with the element
$(0,\dots,0,G^{(n)},0,0,\dots)\in\FF$.

Next, we define the space $\ddG$ of multiple finite
configurations over $X$:
$$\ddG:=\bigsqcup_{n\in\N_0}\overset{{.}{.}}\Gamma{}_X^{(n)}.$$
Here,
$\N_0:=\{0,1,2,\dots\}$,
$\overset{{.}{.}}
\Gamma{}_X^{(0)}:=\{\varnothing\}$ and $\overset{{.}{.}}
\Gamma{}^{(n)}_X$,
$n\in\N$, is the factor space
$$\overset{{.}{.}}\Gamma{}_X^{(n)}:=X^n/S_n$$
with $S_n$ being the group of all permutations of $\{1,\dots,n\}$,
which naturally acts on $X^n$:
\begin{equation}\label{reem}
\sigma(x_1,\dots,x_n)=(x_{\sigma(1)},\dots,x_{\sigma(n)}),\qquad
\sigma\in S_n.\end{equation}
Thus, an  $\eta=[x_1,\dots,x_n]\in \overset{{.}{.}}\Gamma{}_X^{(n)}$
is an equivalence class
 consisting of $n$ elements each of which is a point in $X$
(an $n$-point configuration in $X$ with possibly multiple
points).

Each $\overset{{.}{.}}\Gamma{}_X^{(n)}$ is equipped with the factor
topology generated by the topology on $X^n$, and $\ddG$ is equipped
then by the topology of disjoint union. It follows directly from the
construction of $\FF$ that each $G\in\FF$ can be considered as
the function on $\ddG$ defined by
\begin{align}
G(\varnothing):&=G^{(0)},\notag\\
G([x_1,\dots,x_n]):&=G^{(n)}(x_1,\dots,x_n),\qquad n\in\N.\label{2star}
\end{align}

Notice that in the formula \eqref{2star} we fixed, in fact, a numeration
of the points in $X$ defining the equivalence class, but the
right hand side of \eqref{2star} is independent of the numeration. Now, we will
need a numeration once more to define the notion of summation over
 partitions of an equivalence class.

So, let $\eta=[x_1,\dots,
x_n]$ be an equivalence class with a fixed numeration of points. To each nonempty
subset
$\xi$ of the set $\{1,\dots,n\}$ there corresponds the equivalence class
defined by the points $x_i$, $i\in\xi$.
The $\xi=\emptyset$ as a subset of $\{1,\dots,n\}$ corresponds to the
$\emptyset$ as an element of $\overset{{.}{.}}
\Gamma{}_X^{(0)}$.
 Thus, we will preserve  the
notation $\xi$ for the corresponding element of $\ddG$.
For a function $F\colon (\ddG)^k\to\C$, we let
\begin{equation}\label{2starstar}
\sum_{(\xi_1,\dots,\xi_k)\in{\cal P}_k(\eta)}F(\xi_1,\dots,\xi_k)
\end{equation}
denote the summation over all partitions $(\xi_1,\dots,\xi_k)$ of $\{1,\dots,n\}$.
As easily seen, the result of the summation \eqref{2starstar} is independent
of the numeration.

Now, we define a convolution $\star$ as the mapping
\begin{equation}\label{1}
\star\colon\FF\oplus\FF\to\FF\end{equation}
given by
\begin{equation}\label{2}
(G_1\star G_2)(\eta):=\sum_{(\xi_1,\xi_2,\xi_3)\in{\cal P}_3(\eta)}
G_1(\xi_1\cup\xi_2)G_2(\xi_2\cup\xi_3),
\end{equation}
where ${\cal P}_3(\eta)$ denotes the set of all partitions
$(\xi_1,\xi_2,\xi_3)$ of $\eta$ in 3 parts.

\begin{lem}
\label{lem1}
$\FF$ with the operation $\star$ is a commutative
nuclear algebra.\end{lem}

\noindent{\it Proof}.
For a class $\eta=[x_1,\dots,x_n]$, let $|\eta|:=n$
 Since
$G_1,G_2\in\FF$, there exist  $n_1,n_2\in\N_0$ such that
$G_i(\eta)=0$ if $|\eta|>n_i$. Then, \eqref{2} implies that
$(G_1\star G_2)(\eta)=0$ if $|\eta|>n_1+n_2$.

Next, we note that, for arbitrary $G_1^{(n_1+n_2)}\in\D^{\ho
(n_1+n_2)}$    and $G_2^{(n_2+n_3)}\in\D^{\ho
(n_2+n_3)}$,
the function
\begin{multline*}
G_1^{(n_1+n_2)}(x_1,\dots,x_{n_1},x_{n_1+1},\dots,x_{n_1+n_2})\times\\
\times
G_2^{(n_2+n_3)}(x_{n_1+1},\dots,x_{n_1+n_2},x_{n_1+n_2+1},\dots,
x_{n_1+n_2+n_3})\end{multline*}
belongs to $\D^{\otimes(n_1+n_2+n_3)}$---the $(n_1+n_2+n_3)$-th
tensor power of $\D$---and moreover, it depends continuously on
$G_1$ and $G_2$. Hence, it is easy to see that $G_1\star G_2$
indeed belongs to $\FF$ and that the operation \eqref{1} is continuous.

The commutativity of  $\star$ follows directly from the definition.
Thus, it remains only to show the associativity.
It follows from \eqref2 and an easy combinatoric consideration that
\begin{gather*}
((G_1\star G_2)\star G_3)(\eta)=\sum_{(\xi_1,\xi_2,\xi_3)\in{\cal P}_3(\eta)}
(G_1\star G_2)(\xi_1\cup\xi_2)G_3(\xi_2\cup\xi_3)\\
=\sum_{(\xi_1,\xi_2,\xi_3)\in{\cal P}_3(\eta)}\,\sum_{(\psi_1,\psi_2,\psi_3)
\in{\cal P}_3(\xi_1\cup\xi_2)}
G_1(\psi_1\cup\psi_2)G_2(\psi_2\cup\psi_3)
G_3(\xi_2\cup\xi_3)\\
=\sum_{(\xi_1,\xi_2,\xi_3)\in{\cal P}_3(\eta)}\,
\sum_{(\psi_{11},\psi_{12},\psi_{13})\in{\cal P}_3(\xi_1)}\,
\sum_{(\psi_{21},\psi_{22},\psi_{23})\in{\cal P}_3(\xi_2)}
G_1(\psi_{11}\cup\psi_{12}\cup\psi_{21}\cup\psi_{22})\times\\
\times
G_2(\psi_{12}\cup\psi_{13}\cup\psi_{22}\cup\psi_{23})
G_3(\psi_{21}\cup\psi_{22}\cup\psi_{23}\cup\xi_3)\\
=\sum_{(\xi_1,\dots,\xi_7)\in{\cal P}_7(\eta)}
G_1(\xi_1\cup\xi_4\cup\xi_6\cup\xi_7)G_2(\xi_2\cup\xi_4\cup\xi_5
\cup\xi_7)G_3(\xi_3\cup\xi_5\cup\xi_6\cup\xi_7).
\end{gather*}
Absolutely analogously, one arrives at the same result
when calculating $(G_1\star(G_2\star G_3))(\eta)$.\quad $\blacksquare$\vspace{2mm}

We will need now also the (usual) space of finite configurations
over $X$---denoted by $\Gamman$---which is defined
as a subset of $\ddG$ consisting of $\emptyset$ and
all $\eta=[x_1,\dots,x_n]
\in\ddG$ such that $x_i\ne x_j$ if $i\ne j$.
Each $\eta=[x_1,\dots,x_n]
\in\Gamma_{X,0}$ can be identified with the set $\{x_1,\dots,x_n\}$.
We have $\Gamma_{X,0}=\bigsqcup_{n\in\N_0}\Gamma_X^{(n)}$, where
$\Gamma_X^{(n)}$ is the space of $n$-point configurations in $X$.

The space $\Gamman$
is endowed with the relative topology as a subset of $\ddG$.

Let $\rho$ be a measure on the Borel $\sigma$-algebra ${\cal B}(
\Gamman)$.
Of course, $\rho$ can be considered as a measure on ${\cal B}(\ddG)$
such that the (measurable) set $\ddG\setminus\Gamman$ is of zero $\rho$ measure.
One sees that the restriction of $\rho$ to $\overset{{.}{.}}\Gamma{}
_X^{(n)}$ is actually a measure  on
${\cal B}_{\mathrm sym}(X^n)$. Here,  ${\cal B}_{\mathrm sym}(X^n)$
denotes that sub-$\sigma$-algebra of the Borel $\sigma$-algebra ${\cal B}(X^n)$
consisting of symmetric sets, i.e., sets in $X^n$
which are invariant with respect to the
action \eqref{reem} of the permutation group $S_n$ on $X^n$.
For example, for each Borel $\Lambda\in{\cal B}(X)$ we have
$\Lambda^n\in
{\cal B}_{\mathrm sym}(X^n)$.

We will suppose that $\rho$ satisfies the following
assumptions:

\begin{description}

\item[{\rm (A1)}] {\it Normalization}: $\rho(\Gamma_X^{(0)})=1$.

\item[{\rm (A2)}]
{\it Local finiteness}:
 For each $n\in\N$ and each compact subset $\Lambda
\subset X$, we have $$\rho(\Gamma_{\Lambda}^{(n)})<\infty$$
(where $\Gamma_{\Lambda}^{(n)}$ denotes, of course, the $n$-point
configuration space over $\Lambda$).

\item[{\rm (A3)}] {\it Positive definiteness}: For each $G\in\FF$
$$\int_{\Gamman}(G\star\overline G)(\eta)\,\rho(d\eta)\ge0,
$$
where $\overline G$ is the complex conjugate of $G$.
\end{description}

Thus, it follows from (A2) and (A3) that
$$\FF\oplus\FF\ni
(G_1,G_2)\mapsto a_\rho(G_1,G_2):=\int_{\Gamman}(G_1\star G_2)(\eta)
\,\rho(d\eta)\in\C $$
is a  bilinear continuous form which is positive definite:
$a_\rho(G,\overline G)\ge0$.
Therefore, by using the general technique, e.g.,  \cite{BK}, Ch.~5, Sect.~5,
subsec.~1, we can construct a nuclear factor-space
\begin{equation}\label{pip}
\widehat{{\cal F}}_{\mathrm fin}(\D):=\FF/\{G': a_\rho(G',\overline
G{}')=0\},\end{equation}
consisting of factor classes
$$\widehat G=\{G'\in\FF: a_\rho(G-G',\overline G-\overline G{}')=0\},$$
and then a Hilbert space ${\cal H}_\rho$ as the closure of
$\widehat{{\cal F}}_{\mathrm fin}(\D)$ with respect to the norm generated
by the scalar product $(\widehat G_1,\widehat G_2)_{{\cal H}_\rho}
:=a_\rho(G_1,\overline G_2)$. Thus, as  a result we get a nuclear space
$\widehat{{\cal F}}_{\mathrm fin}(\D)$ that is topologically, i.e.,
densely and continuously, embedded into the Hilbert space ${\cal H}
_\rho$.

 Now, for each $\fii\in\D$, we define an operator
${\cal A}_\fii$ acting on $\FF$ as
$${\cal A}_\fii G:=\fii\star G,\qquad G\in\FF,$$
and let $A_\fii$ be the operator in
${\cal H}_\rho$ with domain $\operatorname{Dom} A_\fii=
\widehat{{\cal F}}_{\mathrm fin}(\D)$ defined by
\begin{equation}\label{3}
A_\fii\widehat G:=\widehat{{\cal A}_\fii G}=
\widehat{\fii\star G},\qquad G\in\FF.
\end{equation}
By Lemma \ref{lem1},
\begin{align*}
a_\rho({\cal A}_\fii G_1,\overline G_2)&=\int_{\Gamman}((\fii\star
G_1)\star\overline G_2)(\eta)\,\rho(d\eta)\\
&=\int_{\Gamman}(G_1\star(\overline{\fii\star G_2}))(\eta)\,\rho(d\eta)\\
&=a_\rho(G_1,\overline{{\cal A}_\fii G_2}),
\end{align*}
and therefore the definition \eqref{3} makes sense due to
\cite{BK}, Ch.~5, Sect.~5, subsec.~2, which uses essentially the Cauchy--Schwartz
inequality.

We strengthen now the condition (A2) by demanding   the following:

\begin{description}

\item[{\rm (A2${}'$)}] For every compact $\Lambda\subset X$, there exists
a constant $C_\Lambda>0$ such that
\begin{equation}
\label{5} \rho(\Gamma_\Lambda^{(n)})\le C_\Lambda^n\qquad \text{for all }
n\in\N.\end{equation}

\end{description}

\begin{lem}\label{lem2} Let \rom{(A1)}, \rom{(A2${}'$),} and \rom{(A3)} hold\rom. Then the
operators $A_\fii$\rom, $\fii\in\D$\rom, with domain
$\widehat{{\cal F}}_{\mathrm fin}(\D)$ are essentially selfadjoint
in ${\cal H}_\rho$ and their closures\rom,  $A^\sim_\fii$\rom,
constitute a family of commuting selfadjoint operators\rom, where
the commutation is understood in the sense of the resolutions
of the identity\rom.\end{lem}

\noindent {\it Proof}. Let us show that, for any  $G^{(n)}\in\F_n(\D)$,
$\widehat G^{(n)}$
is an
analytical vector of each $A_\fii$, i.e., the series
\begin{equation}\label{51}
\sum_{k=0}^\infty \frac{\|A_\fii^k \widehat G{}^{(n)}\|_{{\cal H}_\rho}}{k!}\,
|z|^k,\qquad z\in\C,\end{equation}
has a positive radius of convergence.
So, let us fix $\fii\in \D$ and $G^{(n)}\in\F_n(\D)$ and let $\Lambda$
be a compact set in $X$ such that $\supp\fii\subset\Lambda$ and
$\supp G^{(n)}\subset\Lambda^n$.

We will say that a measurable function $G$
on $\ddG$ has bounded support if there exists a compact set
$\Lambda\subset X$ and $N\in\N$ such that $\supp G\subset \bigsqcup
_{n=0}^N\overset{{.}{.}}\Gamma{}_\Lambda^{(n)}$.
The space of all bounded measurable
functions with bounded support will be denoted by
$B_{\mathrm bs}(\ddG)$.
Evidently, the formula \eqref2 can be extended to the case where
$G_1,G_2\in B_{\mathrm bs}(\ddG)$.

Set now
$$\tilde\fii(\eta):=\sup_{x\in X}|\fii(x)|\,{\bf 1}_\Lambda(\eta),
\qquad \tilde G^{(n)}(\eta):=\sup_{\eta\in\Gamma_X^{(n)}}|G^{(n)}(\eta)|\,
{\bf 1}_{\Lambda^n}(\eta),
$$
where ${\bf 1}_Y(\cdot)$ denotes the characteristic function of a set $Y$.
Denote by $m$ the volume measure on $X$. Without loss of generality,
we can suppose
that $m(\Lambda)\ge1$. Let $\tilde \rho_\Lambda$ be the measure on
$\ddG$ defined by
$$\tilde\rho_\Lambda\restriction\overset{{.}{.}}
\Gamma{}_X^{(n)}:=C_\Lambda^n\, m^{\otimes n},
$$
where $C_\Lambda$ is the constant
from (A2${}'$) corresponding to the set $\Lambda$. Then, it is easy to see that
\begin{align}
\|A_\fii^k \widehat G{}^ {(n)}\|^2_{{\cal H}_{\rho}}&=
\int_{\Gamman}\big((\fii^{\star k}\star G^{(n)})\star
(\fii^{\star k}\star \overline G{}^{(n)})\big)
(\eta)\,\rho(d\eta)\notag\\
&\le\int_{\ddG} \big((\tilde\fii^{\star k}\star \tilde G^{(n)})^
{\star2}\big)(\eta)\,\tilde\rho_\Lambda(d\eta)\notag\\
&=\int_{\ddG}\big((\tilde G^{(n)})^{\star2}\star\tilde\fii^{\star 2k}\big)
(\eta)\,\tilde\rho_\Lambda(d\eta).\label{6}\end{align}

For any $R^{(n)}$ and $f$ from $B_{\mathrm bs}(\ddG)$
which are only not equal to zero  on
$\overset{{.}{.}}\Gamma{}_{X}^{(n)}$ and
$\overset{{.}{.}}\Gamma{}_{X}^{(1)}$, respectively, we have
\begin{equation}
(R^{(n)}\star f)([x_1,\dots,x_k])
=\begin{cases}\sum\limits_{i=1}^{n+1}f(x_i)R^{(n)}([x_1,\dots,\hat x_i,\dots,x_{n+1}]),&
\text{if }k=n+1,\\
\sum\limits_{i=1}^n f(x_i)R^{(n)}([x_1,\dots,x_n]),&\text{if }k=n,\\
0,&\text{otherwise}.\end{cases}
\label{7}\end{equation}
Here, $\hat x_i$ denotes the absence of $x_i$.
Therefore, if additionally $R^{(n)}\ge0$ and $f\ge0$, then
\begin{equation}\label{jjj}
\int_{\ddG}(R^{(n)}\star f)(\eta)\,\tilde\rho_\Lambda(d\eta)
\le C_{\Lambda,f}(2n+1)\int_{\ddG}R^{(n)}(\eta)\, \tilde\rho_\Lambda
(d\eta),\end{equation}
where
$$C_{\Lambda,f}:=\max\big\{
\operatornamewithlimits{ess\,sup}_{x\in X}f(x),C_\Lambda\int_X f(x)\,m(dx)
\big\},$$
which yields
$$
\int_{\ddG}(R\star f)(\eta)\,\tilde\rho_\Lambda(d\eta)
\le 2C_{\Lambda,f}(n+1)\int_{\ddG}R(\eta)\,\tilde\rho_\Lambda
(d\eta)$$
for each $R\in B_{\mathrm bs}(\ddG)$, $R\ge0$, satisfying
$R\restriction{\overset{{.}{.}}\Gamma}{}_X^{(k)}=0$ if $k> n$.

Hence, by using
\eqref{jjj}, we get
\begin{align}
&\int_{\ddG}\big((\tilde G^{(n)})^{\star2}\star\tilde\fii^{\star2k}\big)(\eta)
\,\tilde\rho_\Lambda(d\eta)\notag\\
&\qquad\le 2C_{\Lambda,\tilde\fii}(2n+2k)\int_{\ddG}
\big((\tilde G^{(n)})^{\star2}\star\tilde\fii^{\star(2k-1)}\big)(\eta)
\,\tilde\rho_\Lambda
(d\eta)\notag\\
&\qquad\le(2C_{\Lambda,\tilde\fii})^2(2n+2k)(2n+2k-1)\int_{\ddG}
\big((\tilde G^{(n)})^{\star2}\star\tilde\fii^{\star(2k-2)}(\eta)\big)\,
\tilde\rho
_\Lambda(d\eta)\notag\\
&\qquad\le\dots\le(2C_{\Lambda,\tilde\fii})^{2k}\,\frac{(2n+2k)!}{(2n)!}
\int_{\ddG}(\tilde G^{(n)})^{\star2}(\eta)\,\tilde\rho_\Lambda(d\eta).
\label{more7}
\end{align}

Thus, \eqref6 and \eqref{more7} give
$$\|A_\fii^k \widehat G{}^{(n)}\|_{{\cal H}_\rho}
\le(2C_{\Lambda,\tilde\fii})^k\big((2n)!\big)^{-1/2}2^{n+k}(n+k)!
\,\|\tilde G^{(n)}\|_{{\cal H}_{\tilde\rho_\Lambda}}.$$
Since
$$\sum_{k=0}^\infty \frac{(4C_{\Lambda,\tilde\fii})^k(n+k)!}
{k!}\,|z|^k<\infty\qquad\text{if }|z|<(4C_{\Lambda,\tilde\fii})^{-1}
,$$
the analyticity of $\widehat G{}^{(n)}$ for $A_\fii$ is proven. By using
 Nelson's analytic vector criterium (e.g., \cite{ReedSimon}, Sect.~X.6, or
\cite{BK}, Ch.~5, Th.~1.7)
we conclude that the operators $A_\fii$
are essentially selfadjoint on $\FFF$.

Next, by Lemma \ref{lem1} and \eqref3, the operators $A_{\fii_1}$ and $A_{\fii_2}$,
$\fii_1,\fii_2\in\D$, commute on $\FFF$. Since the operator
$A_{\fii_2}$ is essentially selfadjoint on $\FFF$, the set
$$(A^\sim_{\fii_2}-z\operatorname{id})\FFF,\qquad z\in\C,\,\Im z\ne0,$$
is dense in ${\cal H}_\rho$. Next, again using Lemma~\ref{lem1}
and \eqref3, we get
$$(A^\sim_{\fii_2}-z\operatorname{id})\FFF
=
(({\cal A}_{\fii_2}-z\operatorname{id})\FF)\,\widehat{}
\subset\FFF.$$
Therefore, the operators $A_{\fii_1}^\sim$, $A_{\fii_2}^\sim$, and
$$A_{\fii_1}^\sim\restriction (A_{\fii_2}^\sim-z\operatorname{id})
\FFF$$
have a total set of analytical vectors. Thus, by
\cite{BK}, Ch.~5, Th.~1.15,
the operators commute in the sense of the resolutions of
the identity.\quad$\blacksquare$\vspace{2mm}

Let $\D'$ denote the dual space of $\D$ and let ${\cal C}_\sigma
(\D')$ be the cylinder $\sigma$-algebra on $\D'$ (see e.g.\
\cite{BK}, Ch.~2, Sect.~1, subsec.~9).

\begin{th}\label{th1} Let a measure $\rho$ on $(\Gamman,
{\cal B}(\Gamman))$
satisfy the assumptions \rom{(A1),} \rom{(A2${}'$),} and \rom{(A3).} Then\rom, there exists
a probability measure $\mu$ on $(\D',{\cal C}_\sigma(\D'))$ and
a unitary isomorphism
$$K\colon {\cal H}_\rho\to L^2(\D',{\cal C}_\sigma(\D'),d\mu):=L^2
(d\mu)$$
such that the image of each operator $A^\sim_\fii$\rom, $\fii\in\D$\rom,
under $K$ is the operator of multiplication by the monomial
$\la\fii,\cdot\ra$ in $L^2(d\mu)$:
\begin{equation}\label{8}
KA_\fii^\sim K^{-1}=\la\fii,\cdot\ra\cdot\,.\end{equation}
The unitary $K$ is defined first on the dense
 set $\FFF$ in ${\cal H}_{\rho}$ by the formula
\begin{equation}
\label{9}
\FFF\ni \widehat G=(\widehat G{}^{(n)})_{n=0}^\infty\mapsto K\widehat
G=(K\widehat G)(\omega)=\sum_{n=0}
^\infty \la G^{(n)},\wom n\ra\end{equation}
\rom(the series in \eqref 9 is actually finite\rom) and then it is extended
by continuity to the whole ${\cal H}_\rho$ space\rom. Here\rom,
$G=(G^{(n)})_{n=0}^\infty\in\FF$ is an arbitrary representative
of $\widehat G\in\FFF$, and for any
 $\om\in\D'$\rom, $\wom n\in \D^{'\,\ho n}$ is the $n$-th
Wick power of $\om$ defined by the recurrence relation
\begin{gather}
\wom 0=1,\quad \wom 1=\om,\notag\\
\la\fii^{\ot(n+1)},\wom{(n+1)}\ra=\frac1{n+1}\big[\la
\fii^{\ot(n+1)},\wom n\ho\om\ra-n\la(\fii^2)\ho\fii^{\ot(n-1)},
\wom n \ra\big],
\label{10}\\
\fii\in\D.\notag
\end{gather}
\end{th}

\begin{rem}\label{rem1}\rom{
Let $\DFF$ stand for the dual of $\FF$. This is the topological
direct product of the dual spaces
$\F_n(\D')=\D_{\C}^{'\,\ho n}$ of $\F_n(\D)$. Thus, an arbitrary element
$R$ of $\DFF$ has the form $R=(R^{(n)})_{n=0}^\infty$ where $R^{(n)}\in
{\cal F}_n(\D')$. Next, it follows from \eqref{pip} that the dual
$\DFFF$ of $\FFF$ can be identified with the factor-space
\begin{multline*}\DFFF=\DFF/\big\{\,
R:\,\ll G,R\gg=0\text{ for each $G\in\FF$}\\
\text{such that $a_\rho(G,\overline G)=0$}
\,\big\}.\end{multline*}
Here, $\ll\cdot,\cdot\gg$ denotes the dual pairing between the spaces
$\FF$ and $\DFF$ (as well as the pairing between the spaces
$\FFF$ and $\DFFF$ below). Thus, each element $R\in\DFF$ is a representative
of some element $\widehat R\in\DFFF$. Define now
$$R(\om):=(\wom n)_{n=0}^\infty\in\DFF,$$
and let $\widehat R(\om)$ be the corresponding element of $\DFFF$. Then,
the formula \eqref9 can be rewritten in the form
\begin{equation}\label{suk}
\FFF\ni \widehat G\mapsto K\widehat G=(K\widehat G)(\om)=\ll\widehat G,
\widehat R(\om)\gg,\end{equation}
and hence \eqref{8} yields
\begin{align*}\ll A_\fii\widehat G,\widehat R(\om)\gg&=(K(A_\fii\widehat G))(\om)
=\la\fii,\om\ra(K\widehat G)(\om)\\
&=\la\fii,\om\ra\ll\widehat G,\widehat R(\om)\gg,\qquad \widehat G\in\FFF.
\end{align*}
So, $\widehat R(\om)$ is a generalized joint eigenvector of the family
$A_\fii^\sim$, $\fii\in\D$, belonging to the joint eigenvalue $\om\in\D'$,
and the unitary $K$, written in the form \eqref{suk}, is the Fourier transform
in generalized joint eigenvectors of this family (see \cite{BK}, Ch.~3, for a
detailed exposition of the general theory).
}\end{rem}

\noindent {\it Proof of Theorem}\/ \ref{th1}. We will use
the standard
technique of
construction of the Fourier transform in generalized joint eigenvectors
of a family of commuting selfadjoint operators \cite{BK,
Ly,BeLi}. In fact, the existence of a measure and  a unitary $K$
satisfying \eqref8 and given by the formula \eqref9 with some
kernels $\wom n\in \D^{'\,\ho n}$
follows from the following lemma.

\begin{lem}\label{lem3} \rom{1)} For each $\fii\in\D$\rom,
$A_\fii$ is a linear continuous operator on $\FFF.$

\rom{2)} For an arbitrary fixed $\widehat G\in\FFF$\rom, the mapping
$$\D\ni\fii\mapsto A_\fii \widehat G\in\FFF$$
is linear and continuous\rom.

\rom{3)} The vacuum $\widehat\Omega=(1,0,0,\dots)\,\widehat{}\in\FFF$ is a strong cyclic
vector of the family $(A^\sim_\fii)_{\fii\in\D}$, i.e., the linear span of the set
$$\{\widehat\Omega\}\cup\{A_{\fii_1}\dotsm A_{\fii_n}\widehat\Omega\mid\fii_i\in
\D,\,i=1,\dots,n,\,n\in\N\}$$
is dense in $\FFF.$
\end{lem}

\noindent{\it Proof of Lemma\/} \ref{lem3}. 1) and 2) Clear by Lemma~1.

3) Denote by $\Omega=(1,0,0,\dots)$ the vacuum in $\FF$. It suffices
to show that the set
$$\{\Omega\}\cup\{{\cal A}_{\fii_1}\dotsm {\cal A}_{\fii_n}\Omega\mid\fii_i\in
\D,\,i=1,\dots,n,\,n\in\N\}$$
is dense in $\FF$.

Because of \eqref7, we have on $\FF$
\begin{equation}\label{11}{\cal A}_\fii={\cal A}_\fii^++{\cal A}_\fii^0,\end{equation}
where ${\cal A}_\fii^+$ is a creation operator:
\begin{equation}\label{12}
{\cal A}_\fii^+\psi^{\ot n}=(n+1)\fii\ho\psi^{\ot n},\end{equation}
and ${\cal A}_\fii^0$ is a neutral operator:
\begin{equation}\label{13}
{\cal A}_\fii^0\psi^{\ot n}=n(\fii\psi)\ho\psi^ {\ot(n-1)}.\end{equation}
Therefore, taking to notice that the ${\cal A}_\fii^+$'s are the usual creation
operators,  the cyclicity of $\Omega$ for the operators
${\cal A}_\fii$
follows  from the proof of Theorem~2.1 in \cite{Ly}, p.~65.\quad $\blacksquare$
\vspace{2mm}

To finish the proof of the theorem, we need only to
show that  \eqref{10} holds. To this end, denote for $G\in\FF$ $KG:=K\widehat G$.
Then, upon
\eqref{8}, \eqref{9}, \eqref{11}--\eqref{13},
$$\la\fii,\cdot\ra K(\fii^{\ot n})=K{\cal A}_\fii\fii^{\ot n}
=(n+1)K(\fii^{\ot(n+1)})+nK((\fii^2)\ho\fii^{\ot(n-1)}),$$
which implies  \eqref{10}.\quad$\blacksquare$\vspace{2mm}

\begin{cor}\label{cor1}
Under the conditions of Theorem~\rom{\ref{th1},} we have for each
$G\in{\cal H}_{\rho}$
$$\int_{\Gamman} G(\eta)\,\rho(d\eta)=\int_{\D'}KG(\om)\,\mu(d\om).
$$
\end{cor}

\noindent{\it Proof}. Since $K$ is unitary, we have, for arbitrary $G_1,G_2
\in{\cal H}_\rho$,
$$\int_{\Gamma_{X,0}}(G_1\star\overline G_2)(\eta)\,\rho(d\eta)=\int_{\D'}
(KG_1)(\om)(\overline{KG_2})(\om)\,\mu(d\om).$$
By setting in this formula $G_1=G$ and $G_2=\widehat\Omega$ and noting that,
from one hand side, the vacuum is the identity element for the $\star$-convolution
and on the other hand $K\widehat\Omega\equiv1$, we get the corollary.
\quad $\blacksquare$

\begin{rem}\label{rem2}\rom{Let us consider the functional
$$L(\fii;\om):=e^{\la\log(1+\fii),\om\ra},$$
which is evidently analytical in $\fii$ in a neighborhood of zero
in $\D_\C$ for each fixed $\om\in\D'$. Then, by differentiating
this functional and by using the recurrence relation \eqref{10},
one can show that $L$ is the generating functional of the Wick
monomials
$\la\fii^{\ot n},\wom n\ra$, i.e.,
$$L(\fii,\om)=\sum_{n=0}^\infty\la\fii^{\ot n},\wom n\ra$$
for $\fii$ from a neighborhood of zero (more exactly, for
$\fii\in
\D_\C$ such that
$\sup_{x\in X}|\fii(x)|<1$).
Notice that the functional $L$ is just the character in the generalized
translation operator
approach to Poisson analysis \cite{new2}.
}\end{rem}

\section{The measure $\rho$ as a correlation measure }

The configuration space $\Gamma_X$ over $X$ is defined as the set of
all locally finite subsets (configurations) in $X$:
$$\Gamma_X:=\{\gamma\subset X\mid |\gamma\cap\Lambda|<\infty
\ \text{for each compact }\Lambda\subset X\}.$$
Here $|A|$ denotes the cardinality of a set $A$. One can identify any
$\gamma\in\Gamma_X$ with the positive Radon measure
$$\sum_{x\in\gamma}\delta_x\in {\cal M}(X),
$$
where $\delta_x$ is the Dirac measure with mass in $x$,
 $\sum_{x\in\varnothing}\delta_x:=$zero measure, and
${\cal M}(X)$  stands for the set of all positive
 Radon  measures
on ${\cal B}(X)$. The space $\Gamma_X$
can be endowed with the relative topology as a subset of the space
${\cal M}(X)$ with the vague topology, i.e., the weakest topology on $\Gamma_X$
such that all maps
$$\Gamma_X\ni\gamma\mapsto\la f,\gamma\ra:=\int_X f(x)\,\gamma(dx)
=\sum_{x\in\gamma}f(x)$$
are continuous.
Here, $f\in C_0(X)$ (:$=$the set of all continuous functions in
$X$ with compact support).
We will denote by ${\cal B}(\Gamma_X)$ the Borel $\sigma$-algebra
on $\Gamma_X$. In fact, $\Gamma_X$ is a measurable subset of $\D'$
and the trace $\sigma$-algebra of ${\cal C}_\sigma(\D')$
on $\Gamma_X$ (i.e., the $\sigma$-algebra on $\Gamma_X$ consisting
of intersections of sets from ${\cal C}_\sigma(\D')$
with $\Gamma_X$) coincides with ${\cal B}(\Gamma_X)$.

The following lemma gives a direct representation of the Wick powers $\wom{n}$
in the case where $\omega=\gamma$ is a configuration.

\begin{lem}\label{lem4} For each $\gamma\in\Gamma_X$\rom, we have
\begin{equation}\label{14}
\gom n=\sum_{\eta\Subset\gamma,\,|\eta|=n}\underset{x\in\eta}{\ho}
\delta_x,\end{equation}
where the summation is extended over all $n$-point
subconfigurations from
$\gamma.$\end{lem}

\noindent {\it Proof}.
 For
$n=0$ and $n=1$ the formula evidently holds, and let us suppose that
it holds for all $m\le n$. Then, upon \eqref{10}
\begin{gather*}
\la\fii^{\ot(n+1)},\gom {(n+1)}\ra=\frac1{n+1}\big[
\la\fii^{\ot n},\gom n\ra\la\fii,\gamma\ra-n\la(\fii^2)\ho
\fii^{\ot(n-1)},\gom n\ra
\big] \\
=\frac1{n+1}\bigg(
\sum_{\eta\Subset\gamma,\,|\eta|=n}\,\prod_{y\in\gamma}\fii(y)
-\sum_{\eta\Subset
\gamma,\,|\eta|=n}\,\sum_{x\in\eta}\fii^2(x)\prod_{y\in\eta\setminus\{x\}}
\fii(y)
\bigg)\\
=\frac1{n+1}\sum_{\eta\Subset\gamma,\,|\eta|=n}\,\prod_{x\in\gamma}\fii(x)
\sum_{y\in\gamma\setminus\eta}\fii(y)=\sum_{\eta\Subset
\gamma,\,|\eta|=n}\fii(y).\quad\blacksquare
\end{gather*}\vspace{2mm}

As a direct consequence of Lemma~\ref{lem4} and Corollary~\ref{cor1}, we get

\begin{prop}\label{prop1} Suppose that, under the assumptions of
Theorem~\rom{\ref{th1}}, the measure $\mu$ has the configuration space $\Gamma_X$
as a set of full measure\rom. Then\rom, the operator $K$ coincides with the
$K$-transform between the spaces of functions of
finite and infinite configurations\rom, while
the measure $\rho$ is the correlation measure of $\mu$ \rom{\cite{Len1,
Len2, Len3,Kuna1}.}\end{prop}

To restrict the measure $\mu$ to $\Gamma_X$, we need
an additional condition
on  $\rho$, which is also not very restrictive.

\begin{description}

\item[\rom{(A4)}] Every compact $\Lambda\subset X$ can be covered
by a finite union of open sets $\Lambda_1,\dots,\Lambda_k$, $k\in\N$,
which have compact closures and satisfy the estimate
$$\rho(\Gamma_{\Lambda_i}^{(n)})\le(2+\eps)^{-n}\quad\text{for all
$i=1,\dots,k$ and  $n\in\N$},$$
where $\eps=\eps(\Lambda)>0$.
\end{description}

Suppose, for example, that a measure $\rho$ on $\Gamma_{X,0}$ has density
$\tilde\rho$ with respect to the Lebesgue--Poisson measure
$$\lambda:=\sum_{n=0}^\infty\frac1{n!}\,m^{\otimes n},$$
and suppose that this density fulfills the estimate$$
\operatornamewithlimits{ess\,sup}_{\eta\subset\Gamma_X^{(n)}}\tilde
\rho(\eta)\le n!\,C^n,\qquad
n\in\N,$$
for some constant $C>0$. Then $\rho$ satisfies trivially (A2${}'$) as well as
(A4). (We note that this situation where the measure $\rho$ has density
with respect to the Lebesgue--Poisson measure is typical in applications.)

\begin{th}\label{th2}
Let a measure $\rho$ on $(\Gamman,{\cal B}(\Gamman))$
satisfy the assumptions
\rom{(A1),} \rom{(A2${}'$), (A3), (A4),} and let $\mu$ be the probability measure
on $(\D',
{\cal C}_\sigma(\D'))$ constructed in Theorem~\rom{\ref{th1}.}
Then\rom, $\Gamma_X$ is of full $\mu$
measure\rom.
\end{th}

\noindent {\it Proof}. The proof is a modification of  part of
the proof of Theorem~5.5 in \cite{Kuna1}.

For a function $\fii\in\D_\C$, define a function $e(\fii,\cdot)$
on $\Gamma_{X,0}$ as follows:
$$\Gamma_{X,0}\ni\eta\mapsto e(\fii,\eta):=\prod_{x\in\eta}\fii(x)\in\C.$$
It follows from Remark~\ref{rem2} that
$$e^{\la\fii,\om\ra}=\sum_{n=0}^\infty\la (e^\fii-1)^{\otimes n},\wom n\ra,$$
where $\fii$ belongs to a neighborhood of zero in $\D_\C$, more exactly,
if $\sup_{x\in X}|\fii(x)|<\delta$ for some $\delta>0$. Therefore,
\begin{equation}\label{chych}
e^{\la\fii,\cdot\ra}=Ke(e^\fii-1,\cdot).\end{equation}

Fix a compact $\Lambda\subset X$.
Let ${\cal C}_{\sigma,\Lambda}(\D')$ denote the sub-$\sigma$-algebra
of ${\cal C}_\sigma(\D')$ generated by the functionals of the form
$$\D'\ni\om\mapsto\la\fii,\om\ra\in\C, \qquad
\fii\in\D(\Lambda),$$
where $\D(\Lambda)$ denote the subspace of $\D$ consisting
of those $\fii$ having support in $\Lambda$. Next, let $\mu_\Lambda$
stand for the restriction of the measure $\mu$ to the sub-$\sigma$-algebra
${\cal C}_{\sigma,\Lambda}(\D')$.

Let now $\fii\in\D_\C(\Lambda)$.
It follows from \eqref{chych} that
$$e(e^\fii-1,\cdot)\star e(e^{\overline{\fii}}-1,\cdot)=e(e^{\fii+
\overline\fii}-1,\cdot).$$
Therefore, by using (A2${}'$), we see that there exists $\delta_\Lambda>0
$ such that
$e(e^\fii-1,\cdot)\in{\cal H}_\rho$ provided $\sup_{x\in X}
|\fii(x)|\le\delta_\Lambda$.
Thus, by Corollary~\ref{cor1},
\begin{equation}\label{tfu}
\int_{\D'}e^{\la\fii,\om\ra}\,\mu_\Lambda(d\om)=\int_{\Gamma_{X,0}}
e(e^\fii-1,\eta)\,\rho(d\eta),\qquad \fii\in\D_\C(\Lambda),\
\sup_{x\in X}|\fii(x)|\le\delta_\Lambda.\end{equation}
Thus, the formula \eqref{tfu} gives the analytic extension
of the Fourier transform of the measure
$\mu_\Lambda$ in a neighborhood of zero.

Let us introduce now a mapping ${\cal R}$ which transforms
the set of measurable
functions on $\Gamma_\Lambda$ into itself as follows:
$$({\cal R}F)(\eta):=\sum_{\xi\subset\eta}(-1)^{|\eta\setminus\xi|}
F(\xi),\qquad \eta\in\Gamma_\Lambda.$$
Let now $\Lambda$ satisfy the condition
\begin{equation}\label{chuu}
\rho(\Gamma_\Lambda^{(n)})\le(2+\eps)^n,\qquad \eps>0.\end{equation}
Define on ${\cal B}(\Gamma_\Lambda)$ the set function
$$\widetilde{\mu}_\Lambda(A):=\int_{\Gamma_\Lambda}({\cal R}
{\bf 1}_A)(\eta)\,\rho(d\eta).$$
Since $\sum_{\xi\subset\eta}1=2^n$ if $|\eta|=n$, we conclude that the bound
\eqref{chuu} implies that $\widetilde\mu_{\Lambda}$ is a signed measure.
Therefore, for $\fii\in{\cal D}_\C(\Lambda)$, we have
\begin{equation}\label{huu}
\int_{\Gamma_\Lambda}e^{\la\fii,\eta\ra}\,\widetilde\mu_\Lambda
(d\eta)=\int_{\Gamma_\Lambda}({\cal R}e^{\la\fii,\cdot\ra})(\eta)\,
\rho(d\eta).\end{equation}
Direct calculation shows that
$$({\cal R}e^{\la\fii,\cdot\ra})(\eta)=e(e^\fii-1,\eta),$$
and therefore, we have from \eqref{tfu} and \eqref{huu}
$$\int_{\D'}e^{\la\fii,\om\ra}\,\mu_\Lambda(d\om)=\int_
{\Gamma_\Lambda}e^{\la\fii,\eta\ra}\,\widetilde\mu_\Lambda(d\eta)
,\qquad \fii\in\D_\C(\Lambda),\
\sup_{x\in X}|\fii(x)|\le\delta_\Lambda.
$$
Therefore, $\widetilde\mu_\Lambda$ is a probability measure
on $\Gamma_\Lambda$, and moreover it coincides with the restriction
of the measure $\mu_\Lambda$ to the set $\Gamma_\Lambda$ considered
as a subset of $\D'$.

Hence
\begin{equation}\label{zip}\mu(\widetilde\Gamma_\Lambda)=1,
\end{equation}
where $\widetilde\Gamma_\Lambda$ denotes the set of all $\om\in\D'$
whose restriction to the set $\Lambda$ is a finite sum of delta functions
concentrated in $\Lambda$ and having disjoint support.

Now, let $\Lambda$ be an arbitrary compactum in $X$ and let $\Lambda_1,
\dots,\Lambda_k$ be open subsets of $X$ as in (A4) corresponding to $\Lambda$.
Since
$$\widetilde \Gamma_{\bigcup_{i=1}^k\Lambda_i}=\bigcap_{i=1}^k\widetilde
\Gamma_{\Lambda_i},$$
we conclude that \eqref{zip} holds
 for
{\it each\/} compact $\Lambda\subset X$. From here, we immediately conclude
that $\mu(\Gamma_X)=1$.\quad$\blacksquare$

\begin{center}\bf ACKNOWLEDGMENTS\end{center}

The authors were partially supported by the SFB 256, Bonn
University. Support by the DFG through Projects~436~113/39 and
436~113/43, and by the BMBF through Project~UKR-004-99 is
gratefully
 acknowledged.

\end{document}